\documentclass{article}
\usepackage{amsmath, amsthm, amscd, amsfonts, amssymb, latexsym}
\newtheorem{theorem}{Theorem}[section]      
\newtheorem{lemma}[theorem]{Lemma}         
\newtheorem{cor}[theorem]{Corollary}          
\numberwithin{equation}{section}
\title{The finite time blow-up of the Yang-Mills flow}
\author{Guanxiang Wang and Chuanjing Zhang}

\begin{document}

	\maketitle
	\begin{abstract}
		
		In this paper, we shall prove that, on a non-flat Riemannian vector bundle over a compact Riemannian manifold, the smooth solution of the Yang-Mills flow will blow up in finite time if the energy of the initial connection is small enough. We also consider the finite time blow up for the Yang-Mills flow with the initial curvature near the harmonic form. Furthermore, when $E$ is a holomorphic vector bundle over a compact K\"ahler manifold, then $E$ will admit a projective flat structure if the trace free part of Chern curvature is small enough.
	\end{abstract}
	\section{Introduction}
	Let $E\rightarrow M$ be a vector bundle over a closed $n$-dimensional $(n\geq2)$ Riemannian manifold $(M,g)$ and $H$ be a Riemannian metric on $E$. The Yang-Mills flow was first introduced by Atiyah and Bott\cite{ref1}, it is a time-dependent connection $A=A(t)$ on $E$ solving the following equation
	$$\dfrac{\partial A(t)}{\partial t}=-D_{A(t)}^*F_{A(t)} ,$$
	where $F_{A(t)}$ is the curvature of $A(t)$, $D_{A(t)}^*$ is the adjoint (with respect to a fixed metric) of covariant differential $D_{A(t)}$ on $\mathfrak{g}_E$-valued forms. The Yang-Mills heat flow is the gradient flow of the well-known Yang-Mills functional 
	$$\mathcal{YM}(A)=\int_M|F_A|^2dV_g.$$ 
	The existence and convergence of smooth solutions for the Yang-Mills flow is an essential problem. It was subsequently shown by Daskalopoulos\cite{ref5} for compact Riemannian surface and by Rade\cite{ref18} in dimensions two and three, that the flow exists for all time and converges. In four dimensional case, A.Schlatter\cite{ref19} and M.Struwe\cite{ref20} have studied the global weak solution for the Yang-Mills flow over closed 4-manifolds, not excluding the possibility that point singularities will form within finite time. A.Waldron\cite{ref2} make a progress that in four-dimensional case, the finite-time singularities actually do not occur. He proved that any classical solution of Yang-Mills flow extends smoothly for all time. But the convergence of the flow on four-dimensional manifolds has not been understood well. In the case of holomorphic vector bundle, Donaldson\cite{ref7} proved that the Yang-Mills flow exists smoothly for all time and converges to a Hermitian-Yang-Mills connection when the holomorphic vector bundle is stable.
	\par In general cases or dimensions, the behaviour of the Yang-Mills flow has not been understood well. Feehan\cite{ref17} proved that the existence and convergence of the Yang-Mills flow if the initial value is closed enough to a local minimal connection in some Sobolev space by applying the theory of gradient flow and Lojasiewicz-Simon gradient inequality. Naito\cite{ref11} proved that over $\mathit{S}^d$ $ (d\geq5)$ with its standard round Riemannian metric of radius one, if the bundle is non-flat, then the smooth Yang-Mills flow with small energy initial value will blow up in finite time. Joseph F. Grotowski\cite{ref13} has shown the finite time blow up for a class of $SO(n)$-equivariant initial connections on a trival principal $SO(n)$-bundle over $\mathbb{R}^n$ when dimension $n$ greater than 4. When the Yang-Mills flow admits a long time smooth solution, Hong-Tian\cite{ref9} have analyzed the asymptotic behaviour of the Yang-Mills flow, they showed that the singular set will occur and the set has Hausdorff codimension at least four. A refined structure theorems on the singular set for Yang–Mills flow in dimensions $n\geq4$ is obtaind by Casey Kellehera and Jeffrey Streets\cite{ref3}.  
	\par In this article, we show that if $E$ is a non-flat vector bundle over a closed Riemannian manifold $(M,g)$, then the Yang-Mills flow with small energy initial value will blow up in finite time. This result is a generalization of Naito\cite{ref11} and partially clarifies the behaviour of Yang-Mills flow in higher dimensions. More precisely, we prove the following theorem.
	\begin{theorem}
	Let $(M,g)$ be an $n$-dimensional $(n\geq2)$ smooth closed Riemannian maniflod, $(E,h)$ be a non-flat Riemannian vector bundle over $M$. Then there exists a positive constant $\sigma=\sigma(n,g,E,h)>0$ with the following significance: if $A_0$ is a smooth connection compatible with metric $h$ on $E$, such that $\mathcal{YM}(A_0)<\sigma$, then the smooth solution $A(t)$ for Yang-Mills flow with initial value $A_0$ blows up in finite time.
	\end{theorem}	
        We also consider the finite time blow-up when the curvature is near the harmonic form.  Let $(E,h)$ be a Hermitian vector bundle with rank $r$ over a closed Riemannian manifold $(M,g)$, $\theta$ be the harmonic representation for  $-2\pi c_1(E)$, where $c_1(E)$ is the first Chern class.
	 \begin{theorem}
	 Let $(M,g)$ be an $n$-dimensional $(n\geq2)$ smooth closed Riemannian maniflod, $(E,h)$ be a non-projectively flat Hermitian vector bundle with rank r over $M$. Then there exists a positive constant $\sigma=\sigma(n,g,E,h)>0$ with the following significance: if $A_0$ is a smooth connection compatiable with metric $h$ on $E$, such that
	 	\begin{gather*}
	 		\Vert F_{A_0}-\dfrac{i\theta}{r}\otimes Id_E\Vert_{L^2} <\sigma,
	 	\end{gather*}
	 then the smooth solution $A(t)$ for Yang-Mills flow with initial value $A_0$ blows up in finite time.
	 \end{theorem}
 \par Let $(M,\omega)$ be a compact K\"ahler manifold, $(E, H_0, \overline{\partial}_E)$ be a rank $r$ holomorphic vector bundle with Chern connection $D_{H_0}$. It is well-known that the Yang-Mills flow with initial value $D_{A(0)}=D_{H_0}$ has global solution\cite{ref7}. But, in general, the convergence of the Yang-Mills flow is still unknow unless the holomorphic bundle $E$ is poly-stable. Denote by $\mathcal{A}_{H_0}^{1.1}$ the set of connections compatible with $H_0$ on $E$ and $F_{A_0}^{0.2}=0$. Let
  $$F_{A_0}^\bot=F_{A_0}-\dfrac{1}{r}trF_{A_0}\otimes Id_E,$$
 is the trace free part of the curvature. The following theorem gives a sufficient condition for the convergence of the Yang-Mills flow. This also partially clarifies the asymptotic behaviour of the Yang-Mills flow in holomorphic vector bundle. 
 \begin{theorem}
 	Let $(M,\omega)$ be an $n$-dimensional compact K\"ahler manifold, $(E, K_0,\overline{\partial}_E)$ be a holomorphic vector bundle over $M$. Then there exists a positive constant $\sigma>0$ depending on the geometry of $M$ and $(E, K_0,\overline{\partial}_E)$, with the following significance: if the Chern connection $D_{K_0}$ satisfies $\lVert F_{K_0}^\bot\rVert_{L^2}<\sigma$, then there exists a Hermitian metric $H_0$ which is conformally equivalent to $K_0$, and a Yang-Mills connection $A_\infty\in \mathcal{A}_{H_0}^{1.1}$ satisfying $F_{A_\infty}^\bot=0$.
 \end{theorem}
	\par The paper is organized as follows. In section 2, we will review some basic notations, estimates and some basic results. In section 3, we prove the main results.
	\section{Preliminaries}
	\subsection{Connections and Curvatures on vector bundle}
	 Let us first recall some standard geometric notations and definitions. As before, assume $(M,g)$ is a closed $n$-dimensional Riemannian manifold, $(E,h)$ is a real (or Hermitian) vector bundle with rank $r$ over $M$. A connection $D_A$ on $E$ is a linear differential operator 
	$$D_A:\Gamma(E)\rightarrow\varOmega^1(E)$$
	such that
	$$D_A(f\sigma)=df\otimes\sigma+fD_A\sigma$$
	for all $f\in C^\infty(M)$ and $\sigma\in \Gamma(E)$, where $\Gamma(E)$ is the space of smooth sections of $E$, $\varOmega^p(E)$, is the space of $E$-valued $p$-forms.
	We also require that the connection $A$ is compatible with the metric $h$, i.e.
	$$dh(\gamma,\beta)=h(D_A\gamma,\beta)+h(\gamma,D_A\beta)$$
	for all $\gamma, \beta\in \Gamma(E)$. Suppose $(U_\alpha,\boldsymbol{\varphi}_\alpha)$ is a local trivialization of $E$, the connection takes the form
	$$D_A=d+A_\alpha,$$ 
	where $A_\alpha$ is connection 1-form, it is a matrix valued 1-form. More precisely, since the connection is compatible with the metric on $E$, if we denote that $\mathfrak{g}_E\subseteq End(E)$ is the subbundle of $End(E)$  such that its fibre at $x$ is just the set of skew-symmetric(or skew-Hermitian) endomorphisms of $E_x$ with respect to $h(x)$, then $A_\alpha\in\varOmega^1(\mathfrak{g}_E)$ is $\mathfrak{g}_E$-valued 1-form.
	 The space of connections, which is denoted by $\mathcal{A}_E$, is an affine space
	$$\mathcal{A}_E=D_A+\varOmega^1(\mathfrak{g}_E).$$
	Of course, $D_A$ also induces a connection on $\mathfrak{g}_E$, we also denoted it by $D_A$ for simplicity. Indeed, for any $\phi\in \varOmega^0(\mathfrak{g}_E), \sigma\in\Gamma(E)$, define
	\begin{gather}
		(D_A\phi)\sigma=D_A(\phi(\sigma))-\phi(D_A\sigma).
	\end{gather}
	\par Let $F$ be any vector bundle over $M$, then for each linear connection $D_A$ on $F$, we define an exterior differential 
	$$D_A:\varOmega^p(F)\rightarrow\varOmega^{p+1}(F),$$  
    as follows. For each real valued differential $p$-form $\gamma$, $p\geq 0$, and each smooth section $\sigma$ of $F$, we set
	\begin{align}
		D_A(\gamma\otimes\sigma)=d\gamma\otimes\sigma+(-1)^p\gamma\wedge D_A\sigma,
	\end{align}
	and extend the definition to general $\phi\in \varOmega^p(F)$ by linearity. Combining the connection on $E$ and $\mathfrak{g}_E$ given above, (2.2) gives a exterior differential on $\varOmega^p(E)$ and
	$\varOmega^p(\mathfrak{g}_E)$, $p\geq 0$.
	\par The bundle metric $h$ and Riemannian metric $g$ induce an inner product on $\varOmega^p(E)$. The inner product on the bundle $\varOmega^0(\mathfrak{g}_E)$ is given by the following, for each $a,b\in \varOmega^0(\mathfrak{g}_E)$, 
	$$\langle a,b\rangle:=Tr(ab^{*h}),$$
	where $b^{*h}$ is transpose (or conjugate transpose) of the endomorphism $b$ with respect to the Riemannian (or Hermitian) metric on $E$. Also, the above inner product and $g$ induce an inner product on the bundle $\varOmega^p(\mathfrak{g}_E)$.
	\par For any connection $D_A$ of $E$, its curvature $F_A$ is determined by
	$$F_A=D_A\circ D_A :\Gamma(E)\rightarrow\varOmega^2(E),$$
	which is a $C^\infty$-linear operator on the sections of $E$. More precisely, it is a $\mathfrak{g}_E$-valued 2-form. Locally, the curvature is given by
	\begin{gather}
		F_A=dA_\alpha+A_\alpha\wedge A_\alpha.
	\end{gather}
	The first Bianchi identity $D_AF_A=0$ is familar.
	\subsection{Gauge transformations}
	A gauge transformation $u$ of $E$ is a smooth section of $End(E)$ such that at each $x\in M$, $u(x)$ is an orthogonal or unitary transformation of the fiber $E_x$. The gauge group, denoted by $\mathcal{G}_E$, is the set of gauge transformations.
	There is a natural action of gauge group $\mathcal{G}_E$ on the space of connections $\mathcal{A}_E$: given $u\in \mathcal{G}_E$ and a connection $D_A$, define the action of $u$ on $D_A$ as 
	$$D_{u(A)}=u\circ D_A\circ u^{-1}.$$
	i.e. for each section $\sigma\in \Gamma(E)$, the gauge action is 
	\begin{align}
		D_{u(A)}(\sigma)&=uD_A(u^{-1}(\sigma)).
	\end{align}
	One can easily verify that $D_{u(A)}$ is also a connection on $E$ and its curvature is
	$$F_{u(A)}=uF_Au^{-1}.$$ 
	\subsection{Yang-Mills functional and Yang-Mills flow}
	Given a smooth connection $D_A$ on $E$, we define the energy of the connection $D_A$ by 
	\begin{gather}
		\mathcal{YM}(A)=\int_M|F_A|^2dV_g,
	\end{gather}
	where $dV_g$ is the volume form of Riemannian metric $g$. Its Euler-Lagrange equation is the well-known Yang-Mills equation,
	\begin{equation}
		D_A^*F_A=0.
	\end{equation}
	We call a connection $A$ a Yang-Mills connection if it satisfies the Yang-Mills equation. Since the Yang-Mills functional (2.5) is gauge invariant, the Yang-Mills connection is also gauge invariant.  The Yang-Mills flow with initial value $A_0$ is
	\begin{equation}
		\begin{cases}
			&\dfrac{\partial A(t)}{\partial t}=-D_{A(t)}^*F_{A(t)},\\
			\\
			&A(0)=A_0.
		\end{cases}
	\end{equation}
	It is the $L^2$-gradient flow about the Yang-Mills functional.
	 Let $(E,h)$ be a Hermitian vector bundle with rank $r$ over a closed Riemannian manifold $(M,g)$. Suppose $A$ is a Yang-Mills conection, then by Bianchi identity and (2.6), $trF_A$ is a harmonic form. Assume $\theta$ is the harmonic representation for  $-2\pi c_1(E)$, where $c_1(E)$ is the first Chern class. By Hodge theory, $trF_A=\sqrt{-1}\theta$.
	\subsection{Hermitian-Yang-Mills flow}
	Let $(E,H_0,\overline{\partial}_E)$ be a holomorphic vector bundle over a compact $n$-dimensional K\"ahler manifold $(M,\omega)$, $\mathcal{A}_{H_0}$ be the space of connections compatible with the metric $H_0$ on $E$  and $\mathcal{A}_{H_0}^{1.1}$ be the space of unitary integrable connections of $E$. Denote by $D_{H_0}$ the Chern connection with respect to $H_0$ and $\overline{\partial}_E$. The Hermitian-Yang-Mills flow with initial metric $H_0$ is 
		\begin{equation}
		\begin{cases}
			&H^{-1}(t)\dfrac{\partial H(t)}{\partial t}=-2(\sqrt{-1}\Lambda _\omega F_{H(t)}-\lambda \rm Id_E),\\
			\\
			&H(0)=H_0,
		\end{cases}
	\end{equation}
where $\lambda=\dfrac{2\pi \rm deg(E)}{\rm rank(E)\rm Vol(M,\omega)}$, and  $\rm deg(E)=\int_Mc_1(E)\wedge\dfrac{\omega^{n-1}}{(n-1)!}$. The Hermitian-Yang-Mills(HYM-) flow was first introduced by Donaldson\cite{ref7} and he proved the global existence of the HYM-flow. If the holomorhpic bundle is $\omega$-polystable, the flow converges to the so-called Hermitian-Einstein metric. This is the well-known Donaldson-Uhlenbeck-Yau theorem\cite{ref7}\cite{ref8}\cite{ref14}\cite{ref23}. Let $h=H_0^{-1}H$, then
\begin{equation}
	D_H-D_{H_0}=h^{-1}\partial_{H_0}h,
\end{equation}
\begin{equation}
	F_H-F_{H_0}=\overline{\partial}_E(h^{-1}\partial_{H_0}h),
\end{equation}
and
\begin{equation}
	trF_H=trF_{H_0}+\overline{\partial}\partial\log \det h,
\end{equation}
 where $D_H$ is the Chern connection with respect to $H$ and $\overline{\partial}_E$. Denote the complex gauge group of Hermitian bundle $(E,H_0)$ by $\mathcal{G}^\mathbb{C}$. The group acts on $\mathcal{A}_{H_0}^{1.1}$ as follows: for $\sigma\in \mathcal{G}^\mathbb{C}$,
$$\overline{\partial}_{\sigma(A)}=\sigma\circ\overline{\partial}_A\circ\sigma^{-1},$$
$$\partial_{\sigma(A)}=(\sigma^{*_{H_0}})^{-1}\circ\partial_A\circ\sigma^{*_{H_0}},$$


$$D_{\sigma(A)}=\overline{\partial}_{\sigma(A)}+\partial_{\sigma(A)}.$$


 Choose $\sigma$ such that $\sigma^{*_{H_0}}(t)\circ\sigma(t)=h(t)$. Note that $D_{A_0}=D_{H_0}$. After a direct calculation, one can get
\begin{equation}
       F_{D_{\sigma(A_0)}}=\sigma\circ F_H\circ\sigma^{-1}.
\end{equation}
\par Using the solution $H(t)$ of (2.8), we can construct a solution $A(t)$ for the Yang-Mills flow with initial value $D_{H_0}$\cite{ref7}. In particular, we have


\begin{equation}
      trF_{A(t)}=trF_H=trF_{H_0}+\overline{\partial}\partial\log \det h.
\end{equation}
This is important for us to prove Theorem 1.3.
	\subsection{Basic estimates and results}
	\par The following $\epsilon$-regularity is proved by Chen-Shen\cite{ref4} for Yang-Mills flow, by Hong-Tian\cite{ref9} for Yang-Mills-Higgs flow and Yang-Mills-Higgs case on holomorphic vector bundle over K\"ahler manifolds by Li-Zhang\cite{ref12}, it is crucial for the proof of Theorem 1.1.
\begin{theorem}[$\epsilon$-regularity]
	$\forall\ C>0$, $\exists\  \epsilon_0, \delta_0<1/4$. Assume $A(t)$ is a smooth solution for Yang-Mills flow with initial value $A_0$, and $\mathcal{YM}(A_0)<C$. Then, if for some $0<R<\min\{i_M, \sqrt{t_0}/2\}$, the inequality
	\begin{gather}
		R^{2-n}\int_{P_R(x_0,t_0)}\vert F_A\vert^2dV_gdt<\epsilon_0,
	\end{gather}
	holds, then for any $\delta \in (0,\delta_0)$, we have
	\begin{gather}
		\sup_{P_{\delta R}(x_0,t_0)}\vert F_A\vert^2\leq 16(\delta R)^{-4},
	\end{gather}
	where $P_R(x_0,t_0)=B_R(x_0)\times[t_0-R^2, t_0+R^2]$ and $i_M$ is the infimum of the injectivity radius.
\end{theorem}
In order to analyze the asymptotic behavior of the Yang-Mills flow, we also need the following result\cite{ref9}:
\begin{theorem}[\cite{ref9} Theorem A]
	Let $E$ be a vector bundle over an $n$-dimensional closed Riemannian manifold $M$. Let $A$ be a global smooth solution of Yang-Mills flow in $M\times\left[0,\infty \right)$ with smooth initial value $A_0$. Then there exists a sequence $\{t_i\}$ such that, as $t_i\rightarrow\infty$, $A(x,t_i)$ converges, modulo gauge transformations, to a Yang-Mills connection $A_\infty$ in smooth topology outside a closed set $\varSigma$. And $\mathcal H^{n-4}(\varSigma)$ is finite. Moreover, 
	\begin{gather}
		\varSigma=\bigcap_{\epsilon_0>r>0}\left\{x\in M: \liminf_{k\rightarrow\infty}r^{4-n}\int_{B_r(x)}\vert F_{A(t_k)}\vert^2 dV_g\ge \sigma_1\right\}
	\end{gather}
	for some constants $\epsilon_0, \sigma_1>0$.	
\end{theorem}
In general, the Yang-Mills connection $A_\infty$ and the singular set $\varSigma$ are not unique. If the $C^0$-norm of $F_{A(t)}$ is uniformly bounded along  the flow, then the Yang-Mills connection $A_\infty$ is smooth on the whole manifold. 
\par We recall a basic curvature estimate for Yang-Mills connections, derived by Nakajima\cite{ref10} (also see Tian\cite{ref22} ).
\begin{theorem}[\cite{ref10} Lemma 3.1. \cite{ref22} Theorem 2.21]
	Let $A$ be a Yang-Mills connection of bundle $(E,h)$ over an $n$-dimensional $(n\geq4)$ closed Riemannian manifold M. Then there are $\epsilon=\epsilon(n)>0$ and $C=C(n)>0$, which depend only on $n$ and $M$, such that for any $p\in M$ and $\rho<r_p$, where $r_p$ is a positive constant depending on $p$ and geometry of $M$, whenever
	\begin{gather*}
		\rho^{4-n}\int_{B_\rho(p)}\vert F_A\vert^2dV_g\leq \epsilon,	\end{gather*}
	then
	\begin{gather*}
		\sup_{B_{\rho/4}(p)}\vert F_A\vert \leq \dfrac{C}{\rho^2} \left( \rho^{4-n}\int_{B_\rho(p)} \vert F_A\vert^2 dV_g\right)^{1/2}.
	\end{gather*}
\end{theorem}
Using a finite cover of $M$ by geodesic balls and applying the above theorem, we obtain the following global version.
\begin{cor}
	Let $(E,h)$ be a vector bundle over an $n$-dimensional $(n\geq4)$ closed Riemannian manifold $(M,g)$. Then there exist constants $\epsilon=\epsilon(n)>0$ and $C=C(n)>0$, which depend only on $n$ and $(M,g)$, such that if $A$ is a smooth Yang-Mills connection satisfying
	\begin{gather*}
		\left\|F_A\right\|_{L^2(M)}\leq \epsilon,
	\end{gather*}
	then,
	\begin{gather*}
		\left\|F_A\right\|_{L^\infty(M)}\leq C\left\|F_A\right\|_{L^2(M)}.
	\end{gather*}
\end{cor}
The  $L^{n/2}$-energy gap theorem obtained by Feehan\cite{ref15}\cite{ref16} is also important for the proof of Theorem 1.1.
\begin{theorem}[\cite{ref15}\cite{ref16}. $L^{n/2}$-energy gap]
	Let $(M,g)$ be an $n$-dimensional $(n\ge2)$ closed smooth Riemannian manifold. $(E,h)$ be a real (or Hermitian) vector bundle over $M$. Then there exists a positive constant $\epsilon=\epsilon(n,g,E,h)>0$ with the following significance: if $A$ is a smooth Yang-Mills connection of $E$ with $$\Vert F_A\Vert_{L^{n/2}}<\epsilon,$$ 
	then $A$ is a flat connection.
\end{theorem}
Take together with the Corollary 2.4, we can deduce the following $L^2$-energy gap.
\begin{cor}
Under the same assumption as in Theorem 2.4, then there exists a positive constant $\epsilon=\epsilon(n,g,E,h)>0$ with the following significance: if $A$ is a smooth Yang-Mills connection of $E$ with $\mathcal{YM}(A)<\epsilon$, then $A$ is a flat connection.
\end{cor}	
 Let $(E,h)$ be a Hermitian vector bundle with rank $r$ over a closed Riemannian manifold $(M,g)$, $\theta$ the harmonic representation for  $-2\pi c_1(E)$, where $c_1(E)$ is the first Chern class. With the aid of energy gap for Yang-Mills connections, we can extend the criterion for the existence of flat connections to that for projectively flat connections. Since the principal $PU(r)$-bundle associated to $E$ is flat if and only if the Hermitian bundle $(E,h)$ is projectively flat\cite{ref21}, the proof is trival and we omit it here. More precisely, 
\begin{cor}
	Let $(M,g)$ be an $n$-dimensional $(n\geq2)$ closed smooth Riemannian manifold, $(E,h)$ be a Hermitian vector bundle with rank $r$ over $M$. Then there exists positive a constant $\epsilon=\epsilon(n,g,E,h)>0$, with the following significance. If $A$ is a smooth connection compatible with $h$ on $E$ , and the curvature $F_A$ satisfies
	\begin{gather*}
		D_A^*(F_A-\dfrac{trF_A}{r}\otimes Id_E)=0
	\end{gather*}
	and 
	\begin{gather*}
		\Vert F_A-\dfrac{trF_A}{r}\otimes Id_E \Vert_{L^{n/2}}<\epsilon,
	\end{gather*}
	then $A$ is projectively flat. In particular, if $A$ is a smooth Yang-Mills connection with
	\begin{gather*}
		\Vert F_A-\dfrac{i\theta}{r}\otimes Id_E \Vert_{L^{n/2}}<\epsilon,
	\end{gather*}
	then $F_A=\dfrac{i\theta}{r}\otimes Id_E$.
\end{cor}  
Similar to Corollary 2.6, we can derive the following result:
\begin{cor}
Assume the hypotheses of Theorem 2.3, then there exists a positive constant $\epsilon=\epsilon(n,g,E,h)>0$, with the following significance. If $A$ is a smooth connection compatible with $h$ on $E$, and the curvature $F_A$ satisfies
	\begin{gather*}
		D_A^*(F_A-\dfrac{trF_A}{r}\otimes Id_E)=0
	\end{gather*}
	and 
	\begin{gather*}
		\Vert F_A-\dfrac{trF_A}{r}\otimes Id_E \Vert_{L^2}<\epsilon
	\end{gather*}
	then $A$ is projectively flat. In particular, if $A$ is a smooth Yang-Mills connection with
	\begin{gather*}
		\Vert F_A-\dfrac{i\theta}{r}\otimes Id_E \Vert_{L^2}<\epsilon
	\end{gather*}
	then $F_A=\dfrac{i\theta}{r}\otimes Id_E$.
\end{cor}
The above corollary is important for us to show Theorem 1.2 and Theorem 1.3. To prove Theorem 1.3, we also need the following Lemma,
\begin{lemma}[\cite{ref6}, Lemma 6]
	 Let $(M,\omega)$ be a compact K\"ahler manifold, $(E, K_0, \overline{\partial}_E)$ be a rank $r$ holomorphic vector bundle. Then there exists a Hermitian metric $H_0$ which is conformally equivalent to $K_0$, such that $deth(t)=1$ along the Hermitian-Yang-Mills flow, where $h(t)=H_0^{-1}H(t)$ and $H(t)$ is the solution for Hermitian-Yang-Mills flow with $H(0)=H_0$.
\end{lemma}
Suppose $H_0=e^\phi K_0$, a simply calculation shows that
$$F_{H_0}=F_{K_0}-\partial\overline{\partial}\phi\cdot Id_E,$$
this gives $F_{K_0}^\bot=F_{H_0}^\bot$. Particularly, $\big|F_{K_0}^\bot\big|_{K_0}=\big|F_{H_0}^\bot\big|_{H_0}$ since  $H_0=e^\phi K_0$.
	\section{Proof of the main Results}
	In this section, we will prove Theorem 1.1, 1.2 and 1.3.
	\par Note that, the constant $C$ may be different from line to line.
	\subsection{Proof of Theorem 1.1}
    The main idea of the proof of Theorem 1.1 is to use the $\epsilon$-regularity to deduce the $C^0$-estimate along the Yang-Mills flow. Then by Theorem 2.2 and Corollary 2.6, there must exist a flat connection on $E$. This contradicts that the bundle $E$ is non-flat.
	\begin{proof}[Proof of Theorem 1.1]
		Suppose $A(t)$ is the smooth global solution with initial value $A_0$. By $\epsilon $-regularity, fixing $C=1$, there  exist $ \epsilon_0, \delta_0<1/4$, let $\mathcal{YM}(A_0)<\sigma$ small enough, such that, for large $t_0>0$, the inequality 
		\begin{gather*}
			R_0^{2-n}\int_{P_{R_0}(x,t)}\vert F_A\vert^2dV_gdt<2R_0^{4-n}\sigma<\epsilon_0,
		\end{gather*}
		holds for some $0<R_0<i_M$ and $\forall\ (x,t)\in M\times\left[t_0,\infty\right) $. Then by $\epsilon $-regularity, for any $\delta \in (0,\delta_0)$, 
		\begin{gather*}
			\sup_{P_{\delta R_0}(x,t)}\vert F_A\vert ^2\leq 16(\delta R_0)^{-4}.
		\end{gather*}
		So, along the Yang Mills flow,  $\sup_{M\times \left[0,\infty\right)}\vert F_A\vert^2<C_0<\infty$, for some constant $C_0>0$. According to Theorem 2.2 and (2.16), we know that there exists a sequence $\{t_i\}$ such that, as $t_i\rightarrow \infty$, $A(x,t_i)$ converges, modulo gauge transformations, to a Yang-Mills connection $A$ in smooth topology on the whole manifold $M$ and $\mathcal{YM}(A)<\sigma$, let $\sigma$ small enough, by Corollary 2.6, $A$ is a flat connection, it's impossible since $E$ is non-flat, hence $A(t)$ must blow up in finite time. 
	\end{proof}
	\subsection{Proof of Theorem 1.2}
	The main idea of the proof of Theorem 1.2 is similar to the Theorem 1.1. Denote
	\begin{gather}
		e(A,\theta)=\vert F_A-\dfrac{i\theta}{r}\otimes Id_E\vert ^2.
	\end{gather}
		\begin{proof}[Proof of Theorem 1.2]
		Suppose $A(t)$ is the smooth global solution with initial value $A_0$. Fixing $C_0>0$, let $\sigma$ small enough with $\mathcal{YM}(A_0)<C_0$, there  exist $ \epsilon_0, \delta_0<1/4$. Let $t_0$ large enough and $0<R<i_M$, then for $\forall x_0\in M$, we have
		\begin{align*}
			&R^{2-n}\int_{P_R(x_0,t_0)}|F_{A(t)}|^2dV_gdt\\
			&\leq2R^{2-n}\int_{P_R(x_0,t_0)}\big(e(A,\theta)+\dfrac{1}{r^2}|\theta\otimes Id_E|^2\big)dV_gdt\\
			&\leq2R^{2-n}\int_{t_0-R^2}^{t_0+R^2}\int_Me(A,\theta)dV_gdt+\dfrac{C}{r^2}|\theta\otimes Id_E|_{C^0}^2R^4\\
			&\leq2R^{2-n}\int_{t_0-R^2}^{t_0+R^2}\int_Me(A_0,\theta)dV_gdt+\dfrac{C}{r^2}|\theta\otimes Id_E|_{C^0}^2R^4\\
			&\leq4R^{4-n}\sigma^2+\dfrac{C}{r^2}|\theta\otimes Id_E|_{C^0}^2R^4,
		\end{align*}
	where $C$ is constant depending only on the geometry of $M$. Taking $R$ and $\sigma$ small enough such that $\dfrac{C}{r^2}|\theta\otimes Id_E|_{C^0}^2R^4<\dfrac{\epsilon_0}{2}$ and $4R^{4-n}\sigma^2<\dfrac{\epsilon_0}{2}$, where $\epsilon_0$ is chosen as in the $\epsilon$-regularity. Then by  $\epsilon$-regularity, for any $\delta \in (0,\delta_0)$, 
		\begin{gather*}
			\sup_{P_{\delta R}(x_0,t_0)}|F_{A(t)}|\leq 16(\delta R)^{-4}.
		\end{gather*}
	Since $M$ is compact, we can conclude that there exists a positive constant $C_0$ such that $|F_{A(t)}|_{C^0}\leq C_0<\infty$ for all $t\geq0$. Form Theorem 2.2, we see there exists a sequence $\{t_i\}$ such that, as $t_i\rightarrow \infty$, $A(x,t_i)$ converges, modulo gauge transformations, to a Yang-Mills connection $A_\infty$ in smooth topology on the whole manifold $M$ and $\int_Me(A_\infty,\theta)dV_g<\sigma^2$. Let $\sigma$ be small enough, due to Corollary 2.8, $A$ is a projective-flat connection, which is contradict with that $E$ is non-projective-flat. Hence $A(t)$ must blow up in finite time. 
	\end{proof}
\subsection{Proof of Theorem 1.3}
\begin{proof}[Proof of Theorem 1.3]
	By Lemma 2.9, we can find a Hermitian metric $H_0$, such that
	 $$\big|F_{ H_0}^\bot\big|_{H_0}=\big|F_{K_0}^\bot\big|_{K_0},$$ 
	 and $det h(t)=det H_0^{-1}H(t)=1$ for all $t$, where $F_{H_0}$ is curvature of Chern connection $A_{H_0}$ with respect to $H_0$ and holomorphic structure $\overline{\partial}_E$ and $H(t)$ is long time solution for Hermitian-Yang-Mills flow with $H(0)=H_0$. Let $A(t)$ be the solution for the Yang-Mills flow with $A(0)=A_{H_0}$. Then because of (2.13),
	$$trF_{A(t)}=tr F_{H(t)}=trF_{H_0}+\overline{\partial}\partial\log\det h(t)=trF_{H_0},$$
	and
	\begin{align*}
       \left| F_{A(t)}\right|_{H_0}^2
       &=\left| F_{A(t)}^\bot\right|_{H_0}^2+\dfrac{1}{r^2}\left| trF_{A(t)}\otimes Id_E\right|_{H_0}^2\\
       &=\left| F_{A(t)}^\bot\right|_{H_0}^2+\dfrac{1}{r^2}\left| trF_{H_0}\otimes Id_E\right|_{H_0}^2.
    \end{align*}
Let $t_0$ be large enough, $\forall x_0\in M$ ,  it holds that 
\begin{align*}
&R^{2-2n}\int_{P_R(x_0,t_0)}\left|F_{A(t)}\right|_{H_0}^2dV_gdt\\
&=R^{2-2n}\int_{P_R(x_0,t_0)}\big(\left| F_{A(t)}^\bot\right|_{H_0}^2+\dfrac{1}{r^2}\left| trF_{H_0}\otimes Id_E\right|_{H_0}^2\big)dV_gdt\\
&\leq R^{2-2n}\int_{t_0-R^2}^{t_0+R^2}\int_M\left| F_{A(t)}^\bot\right|_{H_0}^2dV_gdt+R^{2-2n}\int_{P_R(x_0,t_0)}\dfrac{1}{r^2}\left| trF_{H_0}\otimes Id_E\right|_{H_0}^2dV_gdt\\
&\leq R^{2-2n}\int_{t_0-R^2}^{t_0+R^2}\int_M\left| F_{K_0}^\bot\right|_{K_0}^2dV_gdt+\dfrac{C}{r^2}\left| trF_{H_0}\otimes Id_E\right|_{C^0}^2R^4\\
&\leq2R^{4-2n}\sigma^2+\dfrac{C}{r^2}\left| trF_{H_0}\otimes Id_E\right|_{C^0}^2R^4,
\end{align*}
in the second inequality we have used the fact that along the Yang-Mills flow, $\int_M| F_{A(t)}^\bot|_{H_0}^2dV_g$ is non-increasing and $\big|F_{ K_0}^\bot\big|_{K_0}=\big|F_{H_0}^\bot\big|_{H_0}$,
where $C$ is a constant depending only on the geometry of $M$. Take $R$ and $\sigma$ small enough, such that $\dfrac{C}{r^2}\left| trF_{H_0}\otimes Id_E\right|_{C^0}^2R^4\leq \dfrac{\epsilon_0}{2}$ and $2R^{4-2n}\sigma^2\leq  \dfrac{\epsilon_0}{2}$, where $\epsilon_0$ is chosen as in the $\epsilon$-regularity.
Then by $\epsilon$-regularity, for any $\delta \in (0,\delta_0),$
$$
	\sup_{P_{\delta R}(x_0,t_0)}\vert F_A\vert_{H_0}^2\leq 16(\delta R)^{-4}.
$$
Since $M$ is compact, we can conclude that there exists a positive constant $C_0$, such that $\left|F_{A(t)}\right|_{C^0}\leq C_0<\infty$ for all $t\geq0$. Similar to the proof of Theorem 1.1 and 1.2, there exists a sequence $A(t_k)$ converges, modulo gauge transformations, to a Yang-Mills connection $A_\infty\in\mathcal{A}_{H_0}^{1.1}$ in smooth topology on the whole manifold $M$ and $\Vert F_{A_\infty
}^\bot\Vert_{L^2}<\sigma$. Let $\sigma$ small enough, by corollary 2.8, $F_{A_\infty}^\bot=0$. This completes the proof of Theorem 1.3.

\end{proof}

\end{document}